\documentclass[12pt]{article}
\usepackage{amsfonts}
\usepackage{amsfonts,mathrsfs}
\usepackage{CJK,fancyhdr,amscd}
\usepackage{anysize}
\usepackage{graphicx,epsfig}
\usepackage{amsthm,latexsym,amsmath,amssymb}
\usepackage[perpage,symbol*]{footmisc}
\usepackage{cite}

\newtheorem{thm}{Theorem}[section]

\newtheorem{prop}[thm]{Proposition}
\newtheorem{lem}[thm]{Lemma}

\newtheorem{rem}[thm]{Remark}
\def\pf{\noindent{\bf Proof.} }

\makeatletter \@addtoreset{equation}{section} \makeatother
\begin{document}
\baselineskip=20pt  \hoffset=-3cm \voffset=0cm \oddsidemargin=3.2cm
\evensidemargin=3.2cm \thispagestyle{empty}\vspace{10cm}
\parskip 8pt
 \vskip 2mm

\title{\textbf{Infinitely many periodic solutions for second order
Hamiltonian systems} }
\author{ Qingye Zhang $^{\dag}$ $\quad$  Chungen Liu
$^{\ddag}$ \\ {\small \it School of Mathematical Science and LPMC}\\
 {\small\it  Nankai University, Tianjin 300071, People's Republic  of China}}
\date{}
\maketitle

\footnotetext[0]{\footnotesize\;{\it E-mail address}:
zhzhy323@mail.nankai.edu.cn (Q. Zhang), liucg@nankai.edu.cn (C.
Liu).}

\footnotetext[0]{\footnotesize$^\dag$Partially  supported by NFSC of
China(10701043).}

\footnotetext[0]{\footnotesize$^\ddag$Partially  supported by NFSC
of China(10531050,10621101) and 973 Program of STM(2006CB805903).}

\noindent{\bf Abstract} {\small\renewcommand\baselinestretch{1} In
this paper, we study the existence of infinitely many periodic
solutions for second order Hamiltonian systems $\ddot{u}+\nabla_u
V(t,u)=0$, where $V(t, u)$ is either asymptotically quadratic or
superquadratic as $|u|\to \infty$.\par }

\noindent  {\small{\it MSC2000:} 34C25, 37J45, 47J30, 58E05}

\def\<{\langle}
\def\>{\rangle}

\section{Introduction and  main results}
We consider the following second order Hamiltonian systems
\begin{equation}\label{HS}
\left\{\begin{array}{ll}\ddot{u}+\nabla_u V(t,u)=0, &\forall\,t\in\mathbb{R},\\
u(0)=u(T),\, \dot u(0)=\dot u(T),&\;\,T>0,
\end{array}\right.
\end{equation}
where $V\in C^1\left(\mathbb{R}\times \mathbb{R}^N,
\mathbb{R}\right)$  is $T$-periodic and has the form
\begin{equation}\label{Fsplit}
V(t,u)=\frac{1}{2}\<U(t)u,u\>+W(t,u)
\end{equation}
with $U(\cdot)$ a continuous $T$-periodic symmetric matrix. Here
and in the sequel, $\<\cdot,\cdot\>$ and $|\cdot|$ always denote
the standard inner product and the associated norm in
$\mathbb{R}^N$ respectively.

In this paper, we will study the existence of infinitely many
nontrivial solutions of \eqref{HS} via the variant fountain theorems
established in \cite{ZWm2} under the assumption that $W(t,u)$ is
even in $u$, i.e., $W(t,-u)=W(t,u)$ for all
$(t,u)\in[0,T]\times\mathbb{R}^N$. We  divide the problem into the
following two cases.

\subsection{The asymptotically quadratic case}
For the asymptotically quadratic case, we  make the following
assumptions:

\noindent (AQ$_1$) $W(t, u)\ge 0$ for all $(t,u)\in
[0,T]\times\mathbb{R}^N$, and there exist constants $\mu\in (0,2)$
and $R_1>0$ such that
\[
\<\nabla_uW(t,u),u\>\le \mu W(t,u),\quad \forall\, t\in [0,T] \mbox{
and } |u|\ge R_1,
\]
(AQ$_2$) $\lim\limits_{|u|\to 0}\frac{W(t,u)}{|u|^2}=\infty$
uniformly for $t\in [0,T]$, and there exist constants $c_2, R_2>0$
such that
\[
W(t,u)\le c_2|u|,\quad \forall\, t\in [0,T] \mbox{ and } |u|\le R_2,
\]
(AQ$_3$) $\liminf\limits_{|u|\to \infty}\frac{W(t,u)}{|u|}\ge d>0$
uniformly for $t\in [0,T]$.

We state our first main result as follows.
\begin{thm}\label{asymthm}
Assume that {\rm(AQ$_1$)--(AQ$_3$)} hold and that $W(t,u)$ is even
in $u$. Then \eqref{HS} possesses infinitely many nontrivial
solutions.
\end{thm}

\subsection{The superquadratic case}
For the superquadratic case, we  assume

\noindent (SQ$_1$) There exist constants $a_1>0$ and $\nu>2$ such
that
\[
|\nabla_uW(t,u)|\le a_1(1+|u|^{\nu-1}),\quad \forall\, t\in [0,T]
\mbox{ and } u\in\mathbb{R}^N,
\]
(SQ$_2$) $W(t, u)\ge 0$ for all $(t,u)\in [0,T]\times\mathbb{R}^N$,
and $\lim\limits_{|u|\to \infty}\frac{W(t,u)}{|u|^2}=\infty$
uniformly for $t\in [0,T]$,\\
(SQ$_3$) There exist constants $1\le\varrho\in(\nu-2,\infty)$ and
$b>0$ such that
\[
\liminf\limits_{|u|\to
\infty}\frac{\<\nabla_uW(t,u),u\>-2W(t,u)}{|u|^\varrho}\ge b\mbox{
uniformly for } t\in [0,T].
\]

Our second main result reads as follows.
\begin{thm}\label{superthm}
Suppose that {\rm(SQ$_1$)--(SQ$_3$)} are satisfied and that
$W(t,u)$ is even in $u$. Then \eqref{HS} possesses infinitely many
nontrivial solutions.
\end{thm}

With the aid of variational methods, the existence and multiplicity
of periodic solutions for Hamiltonian systems have been extensively
investigated in many papers (see
\cite{AZ,Ant,AM,BB,BC,BR,DL,DG,F,FW,
GM,L,LL,L89,Long,Long90,L90,L95,PT,R78,R80,R83,T98,ZWm1,ZWm3} and
the references therein).

For asymptotically quadratic case, under various twist conditions
via Morse indices or Maslov-type indices, the authors obtained
finitely many periodic solutions in \cite{AZ,LL,Long,ZWm3} without
any evenness assumption, while in the presence of evenness, the
authors in \cite{ZWm1} studied the existence of infinitely many
solutions for \eqref{HS} under the conditions that $W(t,u)$ is
sign-changing and in some sense of at most linear  growth near
infinity, which are totally different from our conditions (AQ$_1$)
and (AQ$_3$) in Theorem \ref{asymthm}. In \cite{DL}, the authors
also obtained infinitely many periodic solutions for first order
Hamiltonian systems. We note that some conditions of Theorem 1.1 in
\cite{DL} will not be satisfied when problem \eqref{HS} with
conditions (AQ$_1$)--(AQ$_3$) is transformed to the corresponding
first order Hamiltonian system  in \cite{DL}.

For the superquadratic case, most of the results on the multiplicity
of periodic solutions were obtained under the so-called
Ambrosetti-Rabinowitz superquadratic condition near infinity with or
without the evenness assumption (see e.g.
\cite{AM,BB,BC,DL,GM,L,L89,Long90,L90,PT,R83}).  As mentioned in
\cite{FW}, for first order Hamiltonian systems, the
Ambrosetti-Rabinowitz superquadratic condition requires the
Hamiltonian $H(t,z)$ to be superquadratic in all components of the
variable $z= (p, q)$, which excludes the case for the second order
Hamiltonian systems \eqref{HS} with $H(t,p, q) = \frac1 2|p|^2 + V
(t,q)$. In \cite{F}, the author introduced a new superquadratic
condition for first order Hamiltonian systems, which requires only a
combined effect of Ambrosetti-Rabinowitz superquadratic nature in
$p$ and $q$ with $z = (p, q)$ and can include the above case for
second order Hamiltonian systems. Subsequently, under the
superquadratic condition of this type, the authors in \cite{FW}
obtained the existence of infinitely many periodic solutions with
the evenness assumption. For second order Hamiltonian systems
\eqref{HS}, we note that the Ambrosetti-Rabinowitz superquadratic
condition is somewhat stronger than the superquadratic condition
given by (SQ$_2$) and (SQ$_3$) in Theorem \ref{superthm}.

\section {Variational setting and proofs of the main results}
In this section, we will first recall some related preliminaries and
establish the variational setting for our problem, and then give the
proofs of the main results.
\subsection{Preliminaries and Variational setting}
Within this subsection, we will introduce the variational setting
for problem \eqref{HS}. Recall that the space $H^1\left(S_T,
\mathbb{R}^N\right)$ becomes a Hilbert space if it is equipped with
the usual norm
\[
\|u\|_1=\left(\int_0^T(|\dot u|^2+|u|^2)dt\right)^{1/2},\quad
\forall\, u \in H^1(S_T, \mathbb{R}^N),
\]
where $S_T=\mathbb{R}/T\mathbb{Z}$.

Denote by $\mathcal {A}$ the operator $-(d^2/dt^2)-U(t)$ on
$L^2\equiv L^2\left((0,T),\mathbb{R}^N\right)$ with domain
$D(\mathcal {A})=H^2\left(S_T, \mathbb{R}^N\right)$. It is known
that $\mathcal {A}$ is a selfadjoint operator with a sequence of
eigenvalues (counted with multiplicity)
\begin{equation}\label{eigenlist}
\lambda_1\leq \lambda_2 \leq \cdots \to\infty
\end{equation} and the
corresponding system of eigenfunctions $\{e_n:
n\in\mathbb{N}\}(\mathcal {A}e_n=\lambda_ne_n)$ forming an
orthogonal basis in $L^2$. Denote by $|\mathcal {A}|$ the absolute
value of $\mathcal {A}$ and let $|\mathcal {A}|^{1/2}$ be the square
root of $|\mathcal {A}|$ with domain $D(|\mathcal {A}|^{1/2})$. By
the elliptic estimate  and Theorem 3.6 in \cite{SR}, we have
\[
D(|\mathcal {A}|^{1/2})=H^1\left(S_T, \mathbb{R}^N\right).
\]
Furthermore, if we define on $H^1(S_T, \mathbb{R}^N)$ a new inner
product and the associated norm by
\begin{align*}
(u, v)_0&=\left(|\mathcal {A}|^{1/2}u, |\mathcal {A}|^{1/2}v\right)_{2}+(u,v)_{2},\\
\|u\|_0&=(u,u)_0^{1/2},
\end{align*}
then $\|\cdot\|_0$ is equivalent to the usual norm $\|\cdot\|_1$ on
$H^1\left(S_T, \mathbb{R}^N\right)$, where $(\cdot, \cdot)_{2}$
denotes the usual inner product on $L^2$. Set
\begin{equation}\label{n-n0}
n^-=\#\{i|\lambda_i<0\},\; n^0=\#\{i|\lambda_i=0\},\;
\bar{n}=n^-+n^0,
\end{equation}
and let
\begin{equation}\label{L-L0L+}
L^2=L^-\oplus L^0\oplus L^+
\end{equation}
be the orthogonal decomposition in $L^2$ with
\begin{equation}\notag
\begin{aligned}
&L^-=\mathrm{span}\{e_1,\ldots,
e_{n^-}\},\;L^0=\mathrm{span}\{e_{n^-+1},\ldots, e_{\bar{n}}\},\\
&L^+=\left(L^-\oplus
L^0\right)^\perp=\overline{\mathrm{span}\{e_{\bar{n}+1},\ldots\}}.
\end{aligned}
\end{equation}
Now we introduce on $H^1\left(S_T, \mathbb{R}^N\right)$ the
following inner product and norm:
\begin{align*}
(u, v)&=\left(|\mathcal {A}|^{1/2}u,|\mathcal {A}|^{1/2}v\right)_{2}+\left(u^0,v^0\right)_{2},\\
\|u\|&=(u,u)^{1/2},
\end{align*}
where $u=u^-+u^0+u^+$ and $v=v^-+v^0+v^+$ with respect to the
decomposition \eqref{L-L0L+}. Let $E=H^1\left(S_T,
\mathbb{R}^N\right)$, then $E$ becomes a Hilbert space with the
inner $(\cdot,\cdot)$. Clearly, norms $\|\cdot\|$ and $\|\cdot\|_0$
are equivalent. Consequently, the norm $\|\cdot\|$ is also
equivalent to the norm $\|\cdot\|_1$ on $E$. From now on
$(E,(\cdot,\cdot),\|\cdot\|)$ becomes our working space.
\begin{rem}\label{orthdecomE}
{\rm It is easy to check that $E$ possesses the orthogonal
decomposition
\begin{equation}\label{Edecom}
E=E^-\oplus E^0\oplus E^+
\end{equation}
with
\begin{equation}\label{E-E+}
E^-=L^-,\;E^0=L^0\quad\mbox{and}\quad E^+=E\cap
L^+=\overline{\mathrm{span}\{e_{\bar{n}+1},\ldots\}}
\end{equation}
where the closure is taken with respect to the norm $\|\cdot\|$.
Evidently, the above decomposition is also orthogonal in $L^2$.}
\end{rem}

By the Sobolev embedding theorem, we get directly the following
lemma.
\begin{lem}\label{Ecptem}
$E$ is compactly embedded in $L^p\equiv L^p\left((0,T),
\mathbb{R}^N\right)$ for $1\leq p\leq \infty$ and hence there exists
$\tau_p>0$ such that
\begin{equation}\label{LpleE}
|u|_p\leq \tau_p\|u\|,\quad\forall\, u\in E,
\end{equation}
where $|\cdot|_p$ denotes the usual norm on $L^p$ for all $1\leq
p\leq \infty$.
\end{lem}

Now we define a functional $\it\Phi$ on $E$ by
\begin{align}\label{funPhi}
{\it\Phi}(u)&=\frac{1}{2}\int_0^T\left(|\dot{u}|^2-\<U(t)u,
u\>\right)dt-{\it\Psi}(u)\notag\\[7pt]
&=\frac{1}{2}\|u^+\|^2-\frac{1}{2}\|u^-\|^2-{\it\Psi}(u)\quad{\rm
where}\;{\it\Psi}(u)=\int_0^TW(t, u)dt
\end{align}
for all $u=u^-+u^0+u^+\in E=E^-\oplus E^0\oplus E^+ $. Note that
(AQ$_1$) and (AQ$_3$) imply
\begin{equation}\label{aWleumu}
W(t,u)\le c_1(1+|u|^\mu),\quad \forall\,(t,u)\in
[0,T]\times\mathbb{R}^N
\end{equation} for some $c_1>0$.
Likewise, by (SQ$_1$), there exists a constant $a_2>0$ such that
\begin{equation}\label{sWleunu}
W(t,u)\le a_1(|u|+|u|^\nu)+a_2,\quad \forall\,(t,u)\in
[0,T]\times\mathbb{R}^N.
\end{equation}
In view of \eqref{aWleumu} (or \eqref{sWleunu}) and Lemma
\ref{Ecptem}, ${\it\Phi}$ and ${\it\Psi}$ are well defined.
Furthermore, we have the following
\begin{prop}\label{propPhi}
Suppose that either {\rm (AQ$_1$)} and {\rm (AQ$_3$)} or {\rm
(SQ$_1$)} is satisfied. Then ${\it\Psi}\in C^1(E, \mathbb{R})$ and
${\it\Psi}':E\to E^*$ is compact, and hence ${\it\Phi}\in C^1(E,
\mathbb{R})$. Moreover,
\begin{align}
{\it\Psi}'(u)v&=\int_0^T\< \nabla_uW(t, u),
v\>dt\label{Psider},\\[5pt]
{\it\Phi}'(u)v&=(u^+, v^+)-(u^-, v^-)-{\it\Psi}'(u)v\notag\\[7pt]
&=(u^+, v^+)-(u^-, v^-)-\int_0^T\< \nabla_uW(t, u),
v\>dt\label{Phider}
\end{align}
for all $u,v\in E=E^-\oplus E^0\oplus E^+ $ with $u=u^-+u^0+u^+$ and
$v=v^-+v^0+v^+$ respectively, and critical points of ${\it\Phi}$ on
$E$ are solutions of \eqref{HS}.
\end{prop}

\pf The proof is standard and we refer to \cite{BR} and \cite{R}.
$\hfill\Box$

Let $E$ be a Banach space with the norm $\|\cdot\|$ and
$E=\overline{\oplus_{j\in \mathbb{N}}X_j}$ with $\dim X_j<\infty$
for any $j\in \mathbb{N}$. Set $Y_k=\oplus_{j=1}^k X_j$ and
$Z_k=\overline{\oplus_{j=k}^\infty X_j}$. Consider the following
$C^1$-functional ${\it\Phi}_\lambda: E\rightarrow \mathbb{R}$
defined by
\[
{\it\Phi}_\lambda(u):=A(u)-\lambda B(u),\; \lambda\in [1,2].
\]
The following two variant fountain theorems  were established in
\cite{ZWm2}.
\begin{thm}[{\cite[Theorem 2.2]{ZWm2}}]\label{avfthm}
Assume that the functional ${\it\Phi}_\lambda$ defined above
satisfies

\noindent {\rm(T$_1$)} ${\it\Phi}_\lambda$ maps bounded sets to
bounded sets uniformly for $\lambda\in [1,2]$, and
${\it\Phi}_\lambda(-u)={\it\Phi}_\lambda(u)$\\$~~~~~~ $ for all
$(\lambda, u)\in [1,2]\times E${\rm ,}

\noindent {\rm(T$_2$)} $B(u)\ge 0$ for all $u\in E$, and $B(u)\to
\infty$ as $\|u\|\to \infty$ on any finite dimensional\\$~~~~~~$
subspace of $E${\rm ,}

\noindent {\rm(T$_3$)} There exist $\rho_k>r_k>0$ such that
\[
\alpha_k(\lambda):=\inf\limits_{u\in
Z_k,\,\|u\|=\rho_k}{\it\Phi}_\lambda(u)\ge
0>\beta_k(\lambda):=\max\limits_{u\in
Y_k,\,\|u\|=r_k}{\it\Phi}_\lambda(u),\quad \forall\,
\lambda\in[1,2]
\]
and
\[
\xi_k(\lambda):=\inf\limits_{u\in
Z_k,\,\|u\|\le\rho_k}{\it\Phi}_\lambda(u)\to 0\quad \mbox{as } k\to
\infty \mbox{ uniformly for } \lambda\in[1,2].
\]

\noindent Then there exist $\lambda_n\to 1$, $u_{\lambda_n}\in
Y_n$ such that
\[
{\it\Phi}'_{\lambda_n}\big|_{Y_n}(u_{\lambda_n})=0,
{\it\Phi}_{\lambda_n}(u_{\lambda_n})\to \eta_k\in
[\xi_k(2),\beta_k(1)]\quad \mbox{as } n\to \infty.
\]
Particularly, if $\{u_{\lambda_n}\}$ has a convergent subsequence
for every $k$, then ${\it\Phi}_1$ has infinitely many nontrivial
critical points $\{u_k\}\subset E\setminus \{0\}$ satisfying
${\it\Phi}_1(u_k)\to 0^-$ as $k\to \infty$.

\end{thm}

\begin{thm}[{\cite[Theorem 2.1]{ZWm2}}]\label{svfthm}
Assume that the functional ${\it\Phi}_\lambda$ defined above
satisfies

\noindent {\rm(F$_1$)} ${\it\Phi}_\lambda$ maps bounded sets to
bounded sets for $\lambda\in [1,2]$, and
${\it\Phi}_\lambda(-u)={\it\Phi}_\lambda(u)$ for all
\\$~~~~~~$ $(\lambda, u)\in [1,2]\times E${\rm ,}

\noindent {\rm(F$_2$)} $B(u)\ge 0$ for all $u\in E${\rm ,} Moreover,
$A(u)\to \infty$ or $B(u)\to \infty$ as $\|u\|\to \infty${\rm ,}

\noindent {\rm(F$_3$)} There exist $r_k>\rho_k>0$ such that
\[
\alpha_k(\lambda):=\inf\limits_{u\in
Z_k,\,\|u\|=\rho_k}{\it\Phi}_\lambda(u)>\beta_k(\lambda):=\max\limits_{u\in
Y_k,\,\|u\|=r_k}{\it\Phi}_\lambda(u),\quad \forall\,
\lambda\in[1,2].
\]
Then
\[\alpha_k(\lambda)\le
\zeta_k(\lambda):=\inf\limits_{\gamma\in \Gamma_k}\max\limits_{u\in
B_k}{\it\Phi}_\lambda(\gamma(u)),\quad \forall\, \lambda\in[1,2],
\]
where $B_k=\{u\in Y_k:\|u\|\le r_k\}$ and $\Gamma_k:=\{\gamma\in
C(B_k, E)\big|\gamma\mbox{ is odd, } \gamma|_{\partial B_k}=id\}$.
Moreover, for a.e. $\lambda\in[1,2]$, there exists a sequence
$\{u_m^k(\lambda)\}_{m=1}^\infty$ such that
\[
\sup\limits_m
\left\|u_m^k(\lambda)\right\|<\infty,\;{\it\Phi}'_\lambda\left(u_m^k(\lambda)\right)\to
0\mbox{ and } {\it\Phi}_\lambda\left(u_m^k(\lambda)\right)\to
\zeta_k(\lambda)\mbox{ as } m\to \infty.
\]
\end{thm}

In order to apply the above two theorems to prove our main
results, we define the functionals $A$, $B$ and
${\it\Phi}_{\lambda}$ on our working space $E$ by
\begin{equation}\label{ABdef}
A(u)=\frac{1}{2}\|u^+\|^2,\quad
B(u)=\frac{1}{2}\|u^-\|^2+\int_0^TW(t, u)dt,
\end{equation}
and
\begin{equation}\label{funlam}
{\it\Phi}_{\lambda}(u)=A(u)-\lambda B(u)
=\frac{1}{2}\|u^+\|^2-\lambda\left(\frac{1}{2}\|u^-\|^2+\int_0^TW(t,
u)dt\right)\\[7pt]
\end{equation}
for all $u=u^-+u^0+u^+\in E=E^-\oplus E^0\oplus E^+ $ and
$\lambda\in [1,2]$. From Proposition \ref{propPhi}, we know that
${\it\Phi}_{\lambda}\in C^1(E, \mathbb{R})$ for all $\lambda\in
[1,2]$.  Let $X_j={\rm span}\{e_j\}$ for all $j\in\mathbb{N}$, where
$\{e_n: n\in\mathbb{N}\}$ is the system of eigenfunctions given
below \eqref{eigenlist}. Note that ${\it\Phi}_1={\it\Phi}$, where
${\it\Phi}$ is the functional defined in \eqref{funPhi}.

\subsection{Proof of Theorem \ref{asymthm}}
In this subsection, we will first establish the following lemmas and
then give a proof of Theorem \ref{asymthm}.
\begin{lem}\label{aBuinfty}
Let {\rm (AQ$_1$)} and {\rm (AQ$_3$)} be satisfied. Then $B(u)\ge 0$
for all $u\in E$ and $B(u)\to \infty$ as $\|u\|\to \infty$ on any
finite dimensional subspace of $E$.
\end{lem}
\pf Evidently, it follows from \eqref{ABdef}  and (AQ$_1$) that
$B(u)\ge 0$ for all $u\in E$.

We claim that for any finite dimensional subspace $F\subset E$,
there exists a constant $\epsilon >0$ such that
\begin{equation}\label{muteps}
m(\{t\in [0,T]:|u(t)|\ge \epsilon \|u\|\})\ge
\epsilon,\quad\forall\, u\in F\setminus \{0\}.
\end{equation}
Here and in the sequel, $m(\cdot)$ always denotes the Lebesgue
measure in $\mathbb{R}$.  If not, for any $n\in \mathbb{N}$, there
exists $u_n\in F\setminus \{0\}$ such that
\[
m\left(\{t\in [0,T]:|u_n(t)|\ge \|u_n\|/n\}\right)< 1/n.
\]
Let $v_n=u_n/\|u_n\|\in F$ for all $n\in\mathbb{N}$. Then
$\|v_n\|=1$ for all $n\in\mathbb{N}$, and
\begin{equation}\label{mv1n}
m(\{t\in [0,T]:|v_n(t)|\ge 1/n\})< 1/n,\quad \forall\,
n\in\mathbb{N}.
\end{equation}
Passing to a subsequence if necessary, we may assume $v_n\to v_0$ in
$E$ for some $v_0\in F$ since $F$ is of finite dimension. Evidently,
$\|v_0\|=1$. In view of Lemma \ref{Ecptem} and the equivalence of
any two norms on $F$ , we have
\begin{equation}\label{vnv0L1}
\int_0^T|v_n-v_0|dt\to 0\quad \mbox{as } n\to \infty
\end{equation}
and
\[
|v_0|_\infty>0.
\]
By the definition of norm $|\cdot|_\infty$, there exists a
constant $\delta_0>0$ such that
\begin{equation}\label{mvo}
m(\{t\in [0,T]:|v_0(t)|\ge \delta_0 \})\ge \delta_0.
\end{equation}
For any $n\in \mathbb{N}$, let
\[
\Lambda_n=\{t\in [0,T]:|v_n(t)|< 1/n\}\quad \mbox{and}\quad
\Lambda_n^c=\mathbb{R}\setminus \Lambda_n=\{t\in [0,T]:|v_n(t)|\ge
1/n\}.
\]
Set $\Lambda_0=\{t\in [0,T]:|v_0(t)|\ge \delta_0 \}$. Then for $n$
large enough, by \eqref{mv1n} and \eqref{mvo}, we have
\begin{align}\label{mLamn0}
m(\Lambda_n\cap\Lambda_0) \ge m(\Lambda_0)-m(\Lambda_n^c) \ge
\delta_0-1/n\ge \delta_0/2.
\end{align}
Consequently, for  $n$ large enough, there holds
\begin{align*}
\int_0^T|v_n-v_0|dt&\ge
\int_{\Lambda_n\cap\Lambda_0}|v_n-v_0|dt\\[5pt]
&\ge \int_{\Lambda_n\cap\Lambda_0}(|v_0|-|v_n|)dt \\[5pt]
&\ge (\delta_0-1/n)\cdot m(\Lambda_n\cap\Lambda_0)\\[3pt]
&\ge \delta_0^2/4>0.
\end{align*}
This is in contradiction to \eqref{vnv0L1}. Therefore \eqref{muteps}
holds.

For the $\epsilon$ given in \eqref{muteps}, let
\[
\Lambda_u=\{t\in [0,T]:|u(t)|\ge \epsilon\|u\|\},\quad\forall\, u\in
F\setminus \{0\}.
\]
Then by \eqref{muteps},
\begin{equation}\label{mLamu}
m(\Lambda_u)\ge \epsilon, \quad\forall\, u\in F\setminus \{0\}.
\end{equation}
By (AQ$_3$), there exists a constant $R_3>R_1$ such that
\begin{equation}\label{halfdR0}
W(t,u)\ge d|u|/2,\quad \forall\,t\in [0,T] \mbox{ and } |u|\ge
R_3,
\end{equation}
where $R_1$ is the constant given in (AQ$_1$). Note that
\begin{equation}\label{utR0}
|u(t)|\ge R_3, \quad\forall\, t\in \Lambda_u
\end{equation}
for any $u\in F$ with $\|u\|\ge R_3/\epsilon$. Combining (AQ$_1$),
\eqref{mLamu} and  \eqref{utR0}, for any $u\in F$ with $\|u\|\ge
R_3/\epsilon$, we have
\begin{align*}
B(u)&=\frac{1}{2}\|u^-\|^2+\int_0^TW(t, u)dt\\[5pt]
&\ge \int_{\Lambda_u}W(t, u)dt\\[5pt]
&\ge \int_{\Lambda_u}d|u|/2 dt\\[5pt]
&\ge d\epsilon\|u\|\cdot m(\Lambda_u)/2\ge d\epsilon^2\|u\|/2.
\end{align*}
This implies $B(u)\to \infty$ as $\|u\|\to \infty$ on any finite
dimensional subspace $F\subset E$. The proof is
complete.$\hfill\Box$

\begin{lem}\label{forT3}
Assume that {\rm (AQ$_1$)}--{\rm (AQ$_3$)} hold. Then there exists a
positive integer $k_1$ and two sequences $0<r_k<\rho_k\to 0$ as
$k\to \infty$ such that
\begin{equation}\label{aalphak}
\alpha_k(\lambda):=\inf\limits_{u\in
Z_k,\,\|u\|=\rho_k}{\it\Phi}_\lambda(u)> 0,\quad \forall\, k\ge
k_1,
\end{equation}
\begin{equation}\label{agammak}
\xi_k(\lambda):=\inf\limits_{u\in
Z_k,\,\|u\|\le\rho_k}{\it\Phi}_\lambda(u)\to 0\quad \mbox{as } k\to
\infty\mbox{ uniformly for } \lambda\in[1,2]
\end{equation}
and
\begin{equation}\label{abetak}
\beta_k(\lambda):=\max\limits_{u\in
Y_k,\,\|u\|=r_k}{\it\Phi}_\lambda(u)<0,\quad \forall\, k\in
\mathbb{N},
\end{equation}
where $Y_k=\oplus_{j=1}^k X_j=\mathrm{span}\{e_1,\ldots, e_k\}$
and $Z_k=\overline{\oplus_{j=k}^\infty
X_j}=\overline{{\mathrm{span}\{e_{k},\ldots\}}}$ for all $k\in
\mathbb{N}$.
\end{lem}
\pf We complete the proof via the following two steps.

\noindent{\it Step 1}. We prove \eqref{aalphak} and \eqref{agammak}.

Note first that $Z_k\subset E^+$ for all $k\ge \bar{n}+1$ by
\eqref{E-E+}, where $\bar{n}$ is the integer defined in
\eqref{n-n0}. By \eqref{LpleE}, for any $u\in E$ with $\|u\|\le
R_2/\tau_\infty$, there holds
\begin{equation}\label{uinfty}
|u|_\infty\le R_2,
\end{equation}
where $R_2$ and $\tau_\infty$ are the constants in (AQ$_2$) and
\eqref{LpleE} respectively. Then for any $k\ge \bar{n}+1$ and
$u\in E^+$ with $\|u\|\le R_2/\tau_\infty$, by (AQ$_2$) and
\eqref{funlam}, we have
\begin{align}\label{aPhilamge}
{\it\Phi}_{\lambda}(u)&\ge\frac{1}{2}\|u\|^2-2\int_0^TW(t,
u)dt\notag\\
&\ge\frac{1}{2}\|u\|^2-2c_2|u|_1, \quad\forall\, \lambda\in [1,2].
\end{align}
Let
\begin{equation}\label{l1}
\ell_k=\sup\limits_{u\in Z_k,\,\|u\|=1}|u|_1,\quad\forall\,
k\in\mathbb{N}.
\end{equation}
Then
\begin{equation}\label{l1to0}
\ell_k\to 0\quad \mbox{as } k\to \infty
\end{equation}
since $E$ is compactly embedded into $L^1$. Consequently,
\eqref{aPhilamge} and \eqref{l1} imply
\begin{equation}\label{gelk}
{\it\Phi}_{\lambda}(u)\ge \frac{1}{2}\|u\|^2-2c_2\ell_k\|u\|
\end{equation}
for all $k\ge \bar{n}+1$ and $u\in E^+$ with $\|u\|\le
R_2/\tau_\infty$. For any $k\in\mathbb{N}$, let
\begin{equation}\label{rhokdef}
\rho_k=8c_2\ell_k.
\end{equation}
Then by \eqref{l1to0}, we have
\begin{equation}\label{rhok0+}
\rho_k\to 0\quad \mbox{as } k\to \infty.
\end{equation}
Evidently, there exists a positive integer $k_1>\bar{n}+1$ such
that
\begin{equation}\label{rhokle1}
\rho_k<R_2/\tau_\infty,\quad \forall\, k\ge k_1.
\end{equation}
For any $k\ge k_1$, \eqref{gelk} together with \eqref{rhokdef} and
\eqref{rhokle1} yields
\[
\alpha_k(\lambda):=\inf\limits_{u\in
Z_k,\,\|u\|=\rho_k}{\it\Phi}_\lambda(u)\ge \rho_k^2/4
>0.
\]
By \eqref{gelk}, for any $k\ge k_1$ and $u\in Z_k$ with $\|u\|\le
\rho_k$, we have
\[
{\it\Phi}_{\lambda}(u)\ge -2c_2\ell_k\rho_k.
\]
Observing that ${\it\Phi}_{\lambda}(0)=0$ by (AQ$_2$), then
\[
0\ge \inf\limits_{u\in
Z_k,\,\|u\|\le\rho_k}{\it\Phi}_\lambda(u)\ge
-2c_2\ell_k\rho_k,\quad \forall\, k\ge k_1.
\]
This together with \eqref{l1to0} and \eqref{rhok0+} implies
\[
\xi_k(\lambda):=\inf\limits_{u\in
Z_k,\,\|u\|\le\rho_k}{\it\Phi}_\lambda(u)\to 0\quad \mbox{as } k\to
\infty \mbox{ uniformly for }
 \lambda\in[1,2].
\]
{\it Step 2}. We show that \eqref{abetak} holds.

For any $k\in \mathbb{N}$, there exists a constant $C_k>0$ such
that
\begin{equation}\label{u1Ck}
|u|_2\ge C_k \|u\|,\quad \forall\,u\in Y_k
\end{equation}
since norms $|\cdot|$ and $\|\cdot\|$ are equivalent on finite
dimensional space $Y_k$. By (AQ$_2$), for any $k\in \mathbb{N}$,
there exists a constant $\delta_k>0$ such that
\begin{equation}\label{Ckdelk}
W(t,u)\ge |u|^2/C_k^2,\quad \forall\, |u|\le \delta_k.
\end{equation}
By \eqref{LpleE}, for any $k\in \mathbb{N}$ and $u\in E$ with
$\|u\|\le \delta_k/\tau_\infty$, it holds
\[
|u|_\infty\le\delta_k,
\]
where $\tau_\infty$ is the constant in \eqref{LpleE}.  Combining
this with \eqref{funlam}, \eqref{u1Ck} and \eqref{Ckdelk}, for any
$k\in \mathbb{N}$ and $u\in Y_k$ with $\|u\|\le
\delta_k/\tau_\infty$, we have
\begin{align}
{\it\Phi}_{\lambda}(u)&\le\frac{1}{2}\|u^+\|^2-\int_0^TW(t,
u)dt\notag\\
&\le\frac{1}{2}\|u\|^2-|u|_2^2\notag/C_k^2\\
&\le\frac{1}{2}\|u\|^2-\|u\|^2=-\frac{1}{2}\|u\|^2\label{Ck},\quad\forall\,\lambda\in
[1,2].
\end{align}
Now for any $k\in \mathbb{N}$, if we choose
\[
0<r_k<\min\left\{\rho_k,\delta_k/\tau_\infty\right\},
\]
then \eqref{Ck} implies
\[
\beta_k(\lambda):=\max\limits_{u\in
Y_k,\,\|u\|=r_k}{\it\Phi}_\lambda(u)\le -r_k^2/2<0,\quad\forall\,
k\in \mathbb{N}.
\]
The proof is complete. $\hfill\Box$

\noindent {\bf Proof of Theorem \ref{asymthm}.} In view of
\eqref{aWleumu}, \eqref{funlam} and Lemma \ref{Ecptem},
${\it\Phi}_\lambda$ maps bounded sets to bounded sets uniformly for
$\lambda\in [1,2]$. Evidently,
${\it\Phi}_\lambda(-u)={\it\Phi}_\lambda(u)$ for all $(\lambda,
u)\in [1,2]\times E$ since $W(t,u)$ is even in $u$. Thus the
condition (T$_1$) of Theorem \ref{avfthm} holds. Lemma
\ref{aBuinfty} shows that the condition (T$_2$) holds, while Lemma
\ref{forT3} implies that the condition (T$_3$) holds for all $k\ge
k_1$, where $k_1$ is given there. Therefore, by Theorem
\ref{avfthm}, for each $k\ge k_1$, there exist $\lambda_n\to 1$,
$u_{\lambda_n}\in Y_n$ such that
\begin{equation}\label{aPhider0}
{\it\Phi}'_{\lambda_n}\big|_{Y_n}(u_{\lambda_n})=0,\;
{\it\Phi}_{\lambda_n}(u_{\lambda_n})\to \eta_k\in
[\xi_k(2),\beta_k(1)]\quad \mbox{as } n\to \infty.
\end{equation}
For the sake of notational simplicity, throughout the remaining
proof of Theorem \ref{asymthm} we always set $u_n=u_{\lambda_n}$ for
all $n\in\mathbb{N}$.

\noindent {\it Claim 1}. $\{u_n\}$ is bounded in $E$.

Indeed, for the constant $R_3$ given in \eqref{halfdR0}, there
exists a constant $M_1>0$ such that
\begin{equation}\label{R3M1}
\left|W(t,u)-\frac{1}{2}\<\nabla_uW(t,u),u\>\right|\le M_1,\quad
\forall\,t\in [0,T] \mbox{ and } |u|\le R_3.
\end{equation}
By virtue of \eqref{Phider}, \eqref{funlam}, \eqref{aPhider0},
\eqref{R3M1} and (AQ$_1$), we have
\begin{align*}
-{\it\Phi}_{\lambda_n}(u_n)&=\frac{1}{2}{\it\Phi}'_{\lambda_n}\big|_{Y_n}(u_n)u_n-
{\it\Phi}_{\lambda_n}(u_n)\\[5pt]
&=\lambda_n\int_0^T\left[W(t,u_n)-\frac{1}{2}\<\nabla_uW(t,u_n),u_n\>\right]dt\\[5pt]
&\ge\lambda_n\int_{\Omega_n}\left[W(t,u_n)-\frac{1}{2}\<\nabla_uW(t,u_n),u_n\>\right]dt-\lambda_n M_1T\\[5pt]
&\ge\frac{\lambda_n(2-\mu)}{2}\int_{\Omega_n}W(t,u_n)dt-\lambda_n M_1T\\[5pt]
&\ge \frac{d\lambda_n(2-\mu)}{4}\int_{\Omega_n}|u_n|dt-\lambda_n
M_1T,\quad \forall\, n\in\mathbb{N},
\end{align*}
where $\Omega_n:=\{t\in [0,T]:|u_n(t)|\ge R_3\}$, and $d,R_3$ are
the constants in \eqref{halfdR0}. Combining this with
\eqref{aPhider0}, there exists a constant $M_2>0$ such that
\begin{equation}\label{unL1M2}
\int_{\Omega_n}|u_n| dt\le M_2,\quad \forall\, n\in\mathbb{N}.
\end{equation}
For any $n\in\mathbb{N}$, let $\chi_n:[0,T]\rightarrow \mathbb{R}$
be the indicator function of $\Omega_n$, that is,
\[
 \chi_n(t) = \begin{cases}
1,&   t\in \Omega_n, \\
0,&  t\notin \Omega_n,
\end{cases}.
\]
Then by the definition of $\Omega_n$ and \eqref{unL1M2}, there hold
\begin{equation}\label{chinun}
\left|(1-\chi_n) u_n\right|_\infty\le R_3\quad \mbox{and}\quad
\left|\chi_n u_n\right|_1\le M_2, \quad \forall\, n\in\mathbb{N}.
\end{equation}
By virtue of  Remark \ref{orthdecomE}, Lemma \ref{Ecptem} and the
H\"{o}lder inequality, we have
\begin{align*}
\left|u_n^-+u_n^0\right|_2^2&=\left(u_n^-+u_n^0,u_n\right)_2\\
&=\left(u_n^-+u_n^0,(1-\chi_n)
u_n\right)_2+\left(u_n^-+u_n^0,\chi_n u_n\right)_2\\
&\le \left|(1-\chi_n) u_n\right|_\infty\left|u_n^-+u_n^0\right|_1
+\left|\chi_n u_n\right|_1
\left|u_n^-+u_n^0\right|_\infty\\
&\le c_3(R_3+M_2)\left|u_n^-+u_n^0\right|_2,\quad \forall\,
n\in\mathbb{N}
\end{align*}
for some $c_3>0$, where the last inequality follows from
\eqref{chinun} and the equivalence of any two norms on finite
dimensional space $E^-\oplus E^0$. Consequently, we get
\[
\left|u_n^-+u_n^0\right|_2\le c_3(R_3+M_2),\quad \forall\,
n\in\mathbb{N}.
\]
In view of the equivalence of norms $\|\cdot\|$ and $|\cdot|_2$ on
$E^-\oplus E^0$ again, there exists a constant $M_3>0$ such that
\begin{equation}\label{u-u0M2}
\left\|u_n^-+u_n^0\right\|\le M_3,\quad \forall\, n\in\mathbb{N}.
\end{equation}
Note that
\[
\left\|u_n^+\right\|^2=2{\it\Phi}_{\lambda_n}(u_n)
+\lambda_n\left\|u_n^-\right\|^2+2\lambda_n\int_0^TW(t,u_n)dt,\quad
\forall\, n\in\mathbb{N}.
\]
Thus by \eqref{LpleE}, \eqref{aWleumu}, \eqref{aPhider0} and
\eqref{u-u0M2}, there holds
\begin{align}\label{aunbound}
\left\|u_n\right\|^2&=\left\|u_n^
+\right\|^2+\left\|u_n^-+u_n^0\right\|^2\notag\\
&=2{\it\Phi}_{\lambda_n}(u_n)
+\lambda_n\left\|u_n^-\right\|^2+\left\|u_n^-+u_n^0\right\|^2+
2\lambda_n\int_0^TW(t,u_n)dt\notag\\
&\le M_4+4c_1|u_n|_\mu^\mu\notag\\
&\le M_4+4c_1\tau_\mu^\mu\|u_n\|^\mu,\quad \forall\,
n\in\mathbb{N}
\end{align}
for some $M_4>0$, where $\tau_\mu$ and and $c_1$ are  the constants
in \eqref{LpleE} and \eqref{aWleumu} respectively. Since $\mu<2$,
\eqref{aunbound} yields that $\{u_n\}$ is bounded in $E$.

\noindent {\it Claim 2}. $\{u_n\}$ possesses a strong convergent
subsequence in $E$.

In fact, by Claim 1, without loss of generality, we may assume
\begin{equation}\label{weaklim}
u_n^-\rightarrow u_0^-,\;u_n^0\rightarrow
u_0^0,\;u_n^+\rightharpoonup u_0^+ \; \mbox{and}\;
u_n\rightharpoonup u_0\quad \mbox{as }n\to \infty
\end{equation}
for some $u_0=u_0^-+u_0^0+u_0^+\in E=E^-\oplus E^0\oplus E^+ $ since
$\dim(E^-\oplus E^0)<\infty$. By virtue of the Riesz Representation
Theorem, ${\it\Phi}'_{\lambda_n}\big|_{Y_n}: Y_n\to Y_n^*$ and
${\it\Psi}':E\to E^*$ can be viewed as
${\it\Phi}'_{\lambda_n}\big|_{Y_n}: Y_n\to Y_n$ and ${\it\Psi}':E\to
E$ respectively, where $Y_n^*$ is the dual space of $Y_n$. Note that
\[
0={\it\Phi}'_{\lambda_n}\big|_{Y_n}(u_n)=
u_n^+-\lambda_n\left(u_n^-+P_n{\it\Psi}'\left(u_n\right)\right),
\quad \forall\, n\in\mathbb{N},
\]
where $P_n:E\rightarrow Y_n$ is the orthogonal projection for all
$n\in\mathbb{N}$, that is,
\begin{equation}\label{unstrlim}
u_n^+=\lambda_n\left(u_n^-+P_n{\it\Psi}'\left(u_n\right)\right),\quad\forall\,
n\in\mathbb{N}.
\end{equation}
By Proposition \ref{propPhi}, ${\it\Psi}':E\to E$ is also compact.
Due to the compactness of ${\it\Psi}'$ and \eqref{weaklim}, the
right-hand side of \eqref{unstrlim} converges strongly in $E$ and
hence $u_n^+\rightarrow u_0^+$ in $E$. Combining this with
\eqref{weaklim}, we have $u_n\rightarrow u_0$ in $E$. Thus Claim 2
is true.

Now from the last assertion of Theorem \ref{avfthm}, we know that
${\it\Phi}={\it\Phi}_1$ has infinitely many nontrivial critical
points. Therefore, \eqref{HS} possesses infinitely many nontrivial
solutions by Proposition \ref{propPhi}. The proof of Theorem
\ref{asymthm} is complete. $\hfill\Box$

\subsection{Proof of Theorem \ref{superthm}}
The following lemmas are needed in the proof of Theorem
\ref{superthm}.
\begin{lem}\label{sBuinfty}
Let {\rm (SQ$_1$)} and {\rm (SQ$_2$)} be satisfied. Then $B(u)\ge 0$
for all $u\in E$. Furthermore, $A(u)\to \infty$ or $B(u)\to \infty$
as $\|u\|\to \infty$.
\end{lem}
\pf Using the similar arguments of the proof of Lemma \ref{aBuinfty}
with (AQ$_1$) and (AQ$_3$) replaced by (SQ$_1$) and (SQ$_2$), we can
prove
\[
B(u)\ge 0,\quad\forall\, u\in E
\]
and
\[
B(u)\to \infty\quad \mbox{as }\|u\|\to \infty
\]
on any finite dimensional subspace of $E$. Consequently,
\[
B(u)\to \infty\quad \mbox{as }\|u\|\to \infty \mbox{ on }
E^-\oplus E^0
\]
since $E^-\oplus E^0$ is of finite dimension. Combining this with
\eqref{Edecom} and \eqref{ABdef}, we have
\[
A(u)\to \infty\mbox{ or } B(u)\to \infty\quad\mbox{as } \|u\|\to
\infty.
\]
The proof is completed. $\hfill\Box$
\begin{lem}\label{forF3}
Assume that {\rm (SQ$_1$)}--{\rm (SQ$_3$)} hold. Then there exists a
positive integer $k_2$ and two sequences $r_k>\rho_k\to \infty$ as
$k\to \infty$ such that
\begin{equation}\label{salphak}
\alpha_k(\lambda):=\inf\limits_{u\in
Z_k,\,\|u\|=\rho_k}{\it\Phi}_\lambda(u)> 0,\quad \forall\, k\ge
k_2
\end{equation}
and
\begin{equation}\label{sbetak}
\beta_k(\lambda):=\max\limits_{u\in
Y_k,\,\|u\|=r_k}{\it\Phi}_\lambda(u)<0,\quad \forall\, k\in
\mathbb{N},
\end{equation}
where $Y_k=\oplus_{j=1}^k X_j=\mathrm{span}\{e_1,\ldots, e_k\}$
and $Z_k=\overline{\oplus_{j=k}^\infty
X_j}=\overline{{\mathrm{span}\{e_{k},\ldots\}}}$ for all $k\in
\mathbb{N}$.
\end{lem}
\pf We divide the proof into two steps.

\noindent{\it Step 1}. We first prove \eqref{salphak}.

By \eqref{sWleunu} and \eqref{funlam}, we have
\begin{align}\label{sPhilamge}
{\it\Phi}_{\lambda}(u)&\ge\frac{1}{2}\|u\|^2-2\int_0^TW(t,
u)dt\notag\\
&\ge\frac{1}{2}\|u\|^2-2a_1\left(|u|_1+|u|_\nu^\nu\right)-2a_2T,\quad
\forall\,(\lambda,u)\in[1,2]\times E^+,
\end{align}
where $a_1,a_2$ are the constants in \eqref{sWleunu}. Let
\begin{equation}\label{lnu}
\ell_\nu(k)=\sup\limits_{u\in Z_k,\,\|u\|=1}|u|_\nu,\quad\forall\,
k\in\mathbb{N}.
\end{equation}
Then
\begin{equation}\label{lnuto0}
\ell_\nu(k)\to 0 \quad \mbox{as } k\to \infty
\end{equation}
since $E$ is compactly embedded into $L^\nu$. Note that
\begin{equation}\label{ZkE+}
Z_k\subset E^+,\quad\forall\,k\ge \bar{n}+1,
\end{equation}
where $\bar{n}$ is the integer given in \eqref{n-n0}. Combining
\eqref{LpleE},  \eqref{sPhilamge}, \eqref{lnu} and \eqref{ZkE+},
for $k\ge \bar{n}+1$, we have
\begin{equation}\label{gel1lnu}
{\it\Phi}_{\lambda}(u)\ge
\frac{1}{2}\|u\|^2-2a_1\tau_1\|u\|-2a_2T-2a_1\ell_\nu^\nu(k)\|u\|^\nu,\quad\forall\,
(\lambda,u)\in [1,2]\times Z_k,
\end{equation}
where $\tau_1$ is the constant given in \eqref{LpleE}. By
\eqref{lnuto0}, there exists a positive integer $k_2\ge\bar{n}+1$
such that
\begin{equation}\label{rkdef}
\rho_k:=(16a_1\ell_\nu^\nu(k))^{1/(2-\nu)}>\max\{16a_1\tau_1+1,16a_2T\},\quad\forall\,k\ge
k_2
\end{equation}
since $\nu<2$. Evidently,
\begin{equation}\label{rhokinfty}
\rho_k\to \infty \quad \mbox{as } k\to \infty.
\end{equation}
Combining \eqref{gel1lnu} and \eqref{rkdef}, direct computation
shows
\[
\alpha_k(\lambda):=\inf\limits_{u\in
Z_k,\,\|u\|=\rho_k}{\it\Phi}_\lambda(u)\ge \rho_k^2/4
>0,\quad\forall\,k\ge
k_2.
\]
{\it Step 2}. We then verify \eqref{sbetak}.

Note that the proof of \eqref{muteps} does not involve the
conditions (AQ$_1$) and (AQ$_3$). Therefore, it still holds here.
Consequently, for any $k\in \mathbb{N}$, there exists a constant
$\epsilon_k>0$ such that
\begin{equation}\label{mLamuk}
m(\Lambda_u^k)\ge \epsilon_k, \quad\forall\,u\in Y_k\setminus \{0\},
\end{equation}
where $\Lambda_u^k:=\{t\in \mathbb{R}:|u(t)|\ge \epsilon_k\|u\|\}$
for all $k\in \mathbb{N}$ and $u\in Y_k\setminus \{0\}$. By {\rm
(SQ$_2$)}, for any $k\in \mathbb{N}$, there exists a constant
$S_k>0$ such that
\begin{equation}\label{Wgeu2}
W(t,u)\ge |u|^2/\epsilon_k^3, \quad\forall\, |u|\ge S_k.
\end{equation}
Combining \eqref{funlam}, \eqref{mLamuk}, \eqref{Wgeu2} and {\rm
(SQ$_2$)}, for any $k\in \mathbb{N}$ and $\lambda\in [1,2]$, we have
\begin{align}
{\it\Phi}_{\lambda}(u)&\le\frac{1}{2}\|u^+\|^2-\int_0^TW(t,
u)dt\notag\\[5pt]
&\le\frac{1}{2}\|u\|^2-\int_{\Lambda_u^k}\left(|u|^2/\epsilon_k^3\right)dt\notag\\[7pt]
&\le\frac{1}{2}\|u\|^2-\epsilon_k^2\|u\|^2m(\Lambda_u^k)/\epsilon_k^3\notag\\[3pt]
&\le\frac{1}{2}\|u\|^2-\|u\|^2=-\frac{1}{2}\|u\|^2\label{Phi-u2}
\end{align}
for all $u=u^-+u^0+u^+\in Y_k$ with $\|u\|\ge S_k/\epsilon_k$. Now
for any $k\in \mathbb{N}$, if we choose
\[
r_k>\max\left\{\rho_k,S_k/\epsilon_k\right\},
\]
then \eqref{Phi-u2} implies
\[
\beta_k(\lambda):=\max\limits_{u\in
Y_k,\,\|u\|=r_k}{\it\Phi}_\lambda(u)\le -r_k^2/2<0,\quad\forall\,
k\in \mathbb{N}.
\]
The proof is complete. $\hfill\Box$

\noindent {\bf Proof of Theorem \ref{superthm}.} It follows from
\eqref{LpleE}, \eqref{sWleunu} and \eqref{funlam} that
${\it\Phi}_\lambda$ maps bounded sets to bounded sets uniformly for
$\lambda\in [1,2]$. In view of the evenness of $W(t,u)$ in $u$, it
holds that ${\it\Phi}_\lambda(-u)={\it\Phi}_\lambda(u)$ for all
$(\lambda, u)\in [1,2]\times E$. Thus the condition (F$_1$) of
Theorem \ref{svfthm} holds. Besides, Lemma \ref{sBuinfty} and Lemma
\ref{forF3} show that the condition (F$_2$) and (F$_3$) hold
respectively for all $k\ge k_2$, where $k_2$ is given in Lemma
\ref{forF3}. Therefore, by Theorem \ref{svfthm}, for any $k\ge k_2$
and a.e. $\lambda\in [1,2]$, there exists a sequence
$\left\{u_m^k(\lambda)\right\}_{m=1}^\infty\subset E$ such that
\begin{equation}\label{aelam}
\sup\limits_m
\|u_m^k(\lambda)\|<\infty,\;{\it\Phi}'_\lambda\left(u_m^k(\lambda)\right)\to
0\mbox{ and } {\it\Phi}_\lambda\left(u_m^k(\lambda)\right)\to
\zeta_k(\lambda)\quad\mbox{as } m\to \infty,
\end{equation}
where
\[
\zeta_k(\lambda):=\inf\limits_{\gamma\in \Gamma_k}\max\limits_{u\in
B_k}{\it\Phi}_\lambda(\gamma(u)),\quad \forall\, \lambda\in[1,2]
\]
with $B_k=\{u\in Y_k:\|u\|\le r_k\}$ and $\Gamma_k:=\{\gamma\in
C(B_k, E)\big|\gamma\mbox{ is odd, } \gamma|_{\partial B_k}=id\}$.
Furthermore, it follows from the proof of Lemma \ref{forF3} that
\begin{equation}\label{zetakb}
\zeta_k(\lambda)\in
\left[\bar{\alpha}_k,\bar{\zeta}_k\right],\quad\forall\,k\ge k_2,
\end{equation}
where $\bar{\zeta}_k:=\max\limits_{u\in B_k}{\it\Phi}_1(u)$ and
$\bar{\alpha}_k:=\rho_k^2/4\to \infty$ as $k\to \infty$ by
\eqref{rhokinfty}. Using the similar arguments of the proof of Claim
2 in the proof of Theorem \ref{asymthm}, for each $k\ge k_2$, we can
choose $\lambda_n\to 1$ such that the sequence
$\left\{u_m^k(\lambda_n)\right\}_{m=1}^\infty$ obtained by
\eqref{aelam} has a strong convergent subsequence. Without loss of
generality, we may assume
\[
\lim\limits_{m\to \infty}u_m^k(\lambda_n)=u_n^k\in E,\quad
\forall\,n\in\mathbb{N} \mbox{ and } k\ge k_2.
\]
This together with \eqref{aelam} and \eqref{zetakb} yields
\begin{equation}\label{sPhider0}
{\it\Phi}'_{\lambda_n}(u_n^k)=0,\; {\it\Phi}_{\lambda_n}(u_n^k)\in
\left[\bar{\alpha}_k,\bar{\zeta}_k\right],\quad \forall\, n\in
\mathbb{N} \mbox{ and } k\ge k_2.
\end{equation}

As in the proof of Theorem \ref{asymthm}, we claim that the sequence
$\{u_n^k\}_{n=1}^\infty$ in \eqref{sPhider0} is bounded in $E$ and
possesses a strong convergent subsequence with the limit $u^k\in E$
for all $k\ge k_2$. In fact, by (SQ$_3$), there exist constants
$L_0>0$ and $D_1>0$ such that
\begin{equation}\label{L3def}
\frac{1}{2}\<\nabla_uW(t,u),u\>-W(t,u)\ge
\frac{b}{4}|u|^\varrho,\quad \forall\, t\in [0,T] \mbox{ and }
|u|\ge L_0
\end{equation}
and
\begin{equation}\label{D1def}
\left|\frac{1}{2}\<\nabla_uW(t,u),u\>-W(t,u)\right|\le D_1,\quad
\forall\, t\in [0,T] \mbox{ and } |u|\le L_0.
\end{equation}
For notational simplicity, we will set $u_n=u_n^k$ for all $n\in
\mathbb{N}$ throughout this paragraph. By virtue of \eqref{Phider},
\eqref{funlam} and \eqref{sPhider0}-- \eqref{D1def}, we have
\begin{align*}
{\it\Phi}_{\lambda_n}(u_n)&={\it\Phi}_{\lambda_n}(u_n)-\frac{1}{2}{\it\Phi}'_{\lambda_n}(u_n)u_n
\\[5pt]
&=\lambda_n\int_0^T\left[\frac{1}{2}\<\nabla_uW(t,u_n),u_n\>-W(t,u_n)\right]dt\\[5pt]
&\ge\lambda_n\int_{\Pi_n}\left[\frac{1}{2}\<\nabla_uW(t,u_n),u_n\>-W(t,u_n)\right]dt-\lambda_nD_1T\\[5pt]
&\ge\frac{b\lambda_n}{4}\int_{\Pi_n}|u_n|^\varrho
dt-\lambda_nD_1T, \quad \forall\, n\in\mathbb{N},
\end{align*}
where $\Pi_n:=\{t\in [0,T]:|u_n(t)|\ge L_0\}$. This together with
\eqref{sPhider0} implies that
\begin{equation}\label{PinD2}
\int_{\Pi_n}|u_n|^\varrho dt\le D_2, \quad \forall\, n\in\mathbb{N}
\end{equation}
for some $D_2>0$. Since $\varrho\ge1$, it also holds
\begin{equation}\label{PinD3}
\int_{\Pi_n}|u_n|dt\le D_3, \quad \forall\, n\in\mathbb{N}
\end{equation}
for some $D_3>0$ by using the H\"{o}lder inequality if necessary.
Then the similar arguments of the proof of Claim 1 in the proof of
Theorem \ref{asymthm} yields
\begin{equation}\label{u-u0D4}
\left\|u_n^-+u_n^0\right\|\le D_4,\quad \forall\, n\in\mathbb{N}
\end{equation}
for some $D_4>0$. Combining \eqref{LpleE}, \eqref{Phider},
\eqref{PinD2}, \eqref{u-u0D4} and  (SQ$_1$), we have
\begin{align}\label{sunbound}
\|u_n\|^2&=\left\|u_n^
+\right\|^2+\left\|u_n^-+u_n^0\right\|^2\notag\\[5pt]
&={\it\Phi}'_{\lambda_n}(u_n)u_n
+\lambda_n\left\|u_n^-\right\|^2+\left\|u_n^-+u_n^0\right\|^2+
\lambda_n\int_0^T \<\nabla_uW(t,u_n),u_n\>dt\notag\\[5pt]
&\le D_5+2\int_0^Ta_1(1+|u_n|^{\nu-1})|u_n|dt\notag\\[5pt]
&= D_5+2a_1|u_n|_1+2a_1\int_0^T|u_n|^\nu dt\notag\\[5pt]
&\le D_5+2a_1\tau_1\|u_n\|+2a_1\left[\int_{[0,T]\setminus
\Pi_n}|u_n|^\nu dt+\int_{\Pi_n}|u_n|^\nu dt\right]\notag\\[5pt]
&\le
D_5+2a_1TL_0^\nu+2a_1\tau_1\|u_n\|+2a_1|u_n|_\infty^{\nu-\varrho}\int_{\Pi_n}|u_n|^\varrho
dt\notag\\[5pt]
&\le
D_5+2a_1TL_0^\nu+2a_1\tau_1\|u_n\|+2a_1D_2\tau_\infty^{\nu-\varrho}\|u_n\|^{\nu-\varrho}
\end{align}
for some $D_5>0$, where $\tau_1$ and $\tau_\infty$ are the constants
in \eqref{LpleE}. Since $\varrho>\nu-2$, \eqref{sunbound} implies
that $\{u_n\}$ is bounded in $E$. By virtue of the similar arguments
of the proof of Claim 2 in the proof of Theorem of \ref{asymthm}
again, we see that $\{u_n\}$ has a strong convergent subsequence.

Now  for each $k\ge k_2$, by \eqref{sPhider0}, the limit $u^k$ is
just a critical point of ${\it\Phi}={\it\Phi}_1$ with
${\it\Phi}(u^k)\in \left[\bar{\alpha}_k,\bar{\zeta}_k\right]$. Since
$\bar{\alpha}_k\to \infty$ as $k\to \infty$ in \eqref{zetakb}, we
get infinitely many nontrivial critical points of ${\it\Phi}$.
Therefore, \eqref{HS} possesses infinitely many nontrivial solutions
by Proposition \ref{propPhi}. The proof of Theorem \ref{asymthm} is
complete. $\hfill\Box$

\end{document}